\documentclass[11pt]{amsart}

\usepackage{amsthm}
\usepackage{amsmath}
\usepackage{amssymb}
\usepackage{color}

\newtheorem{thm}{Theorem}[section]
\newtheorem{theorem}[thm]{Theorem}

\newtheorem{hadamard}[thm]{Hadamard Conjecture}
\newtheorem{obs}[thm]{Observation}

\newtheorem{definition}[thm]{Definition}
\theoremstyle{remark}

\begin{document}

\title{Partitioning Hadamard vectors into hadamard
matrices}
\author[Casazza and Tremain
 ]{Peter G. Casazza and Janet C. Tremain}
\address{Department of Mathematics, University
of Missouri, Columbia, MO 65211-4100}

\thanks{The authors were supported by
 NSF DMS 1307685; NSF ATD 1321779; and ARO  W911NF-16-1-0008}

\email{Casazzap@missouri.edu; Tremainjc@missouri.edu}


\begin{abstract}
We will show that in a space of dimension $m$, any 
family of $2^{m-1}$ distinct Hadamard vectors (where
you can choose x or -x but not both) can be
partitioned into Hadamard matrices if and only if
$m=2^n$ for some n. We will
solve this problem with a simple algorithm for assigning
the vectors to the Hadamard matrices.  
\end{abstract}

\maketitle

\section{Introduction}

Recall that 
\begin{definition}
In $R^n$, a vector of the form $x=(\pm1,\pm1,\ldots,
\pm1)$ is called a {\bf Hadamard vector}.  An orthogonal
matrix made up of Hadamard vectors is called a
{\bf Hadamard matrix}.
\end{definition}

Note that there are $2^n$ Hadamard vectors but they come
as pairs $x,-x$ and so their are, up to sign, $2^{n-1}$ {\bf distinct Hadamard vectors}.  The {\bf Hadamard Conjecture} states:

\begin{hadamard}
There exists a $(4n)\times (4n)$ Hadamard matrix 
for every n.
\end{hadamard}

These are the only possible cases since for $n=2m+1$,
a maximal set of orthogonal Hadamard vectors contains
just two vectors.

In this paper we will prove the following theorem.

\begin{theorem}
In a space of dimension $m=2^n$, any maximal set of distinct Hadamard vectors
can be partitioned into $2^{2^n-n-1}$ Hadamard matrices.
Moreover, this result fails for all other values of m.
\end{theorem}

\section{Some Elemntary Observations}

We make a few simple observations.

\begin{obs}
If we can prove that just one choice of a maximal set
of distinct Hadamard vectors can be partitioned into
Hadamard matrices, then this is also true for any choice
of distinct Hadamard vectors.
\end{obs}
This is clear since any second choice of distinct 
Hadamard vectors has the property that for any vector
x in this set, either x or -x is in the first set.  And
so the partition of the first set is a partition of the
second set with perhaps sign changes of some rows of
the Hadamard matrix - which is still an orthogonal matrix.

\begin{obs}
The moreover part of the theorem is essentially obvious.
\end{obs}
Given any $m=4n$, if we can partition the $2^{4n-1}$ distinct
Hadamard vectors into $(4n)\times (4n)$ Hadamard matrices,
let $k$ be the number of such matrices.  Then $k\cdot (4n)$
uses up all the distinct Hadamard vectors and so 
\[ k\cdot (4n) = 2^{4n-1},\mbox{ and so m divides }
2^{4n-1}.\]

We will do the proof by induction on n with the case
n=2 below:
\[ A= \begin{bmatrix}
+&+&+&+\\+&+&-&-&\\
+&-&+&-&\\
+&-&-&+
\end{bmatrix}\ \ \ \ \ \ \ \
B=\begin{bmatrix}
+&+&+&-\\
+&+&-&+\\
+&-&+&+\\
+&-&-&-&
\end{bmatrix}\]

\section{Proof of the Theorem}

We assume the theorem holds for $m=2^n$ and we have
partitioned a distinct set of Hadamard vectors into Hadamard matrices
$\{A_i\}_{i=1}^{2^{2^n}-n-1}$ and let $\{x_{ij}\}_{j=1}^{2^n}$ be the row vectors of $A_i$.  We will construct
$2^{2^{n+1}-(n+1)-1}$ matrices of distinct Hadamard
vectors of order $2^n \times 2^{n+1}$ so that
cutting each
of these matrices vertically in half, each of the left
halves and the right halves are orthogonal matrices. 
For each of these, say $[A\ \ B]$, we then take:
\begin{equation}\label{E}\begin{bmatrix}
A&B\\
A&-B
\end{bmatrix}
\end{equation}
and have a partition of distinct Hadamard vectors into
Hadamard matrices.  Since the total number of vectors here
is equal to the total number of distinct Hadamard vectors
for $2^{2^{n+1}}$, we are done.  It will be obvious from our
construction that the vectors we construct are unique.

For any $n$, we define the row shift of an $n\times n$
matrix $A$ with row vectors $\{x_i\}_{i=1}^n$ by:
\[ T_n=\begin{bmatrix}
0&0&0&\cdots&0&1\\
1&0&0&\cdots&0&0\\
\vdots&\vdots&\vdots&\vdots&\vdots&\vdots\\
0&0&0&\cdots&1&0
\end{bmatrix}
\mbox{ so that } T_n\begin{bmatrix}
x_1\\
x_2\\
x_3\\
\vdots\\x_n
\end{bmatrix}
= \begin{bmatrix}
x_n\\
x_1\\
x_2\\
\vdots\\
x_{n-1}
\end{bmatrix}\]

For the construction, for each $A_i,A_j$ above, we form
the $2^{2^n} \times 2^{2^n+1}$  matrices:
\[ \begin{bmatrix}
A_i&\vline &T_m^kA_j
\end{bmatrix}
\mbox{ for all }k=1,2,\ldots,2^n.\]
Since the rows of $A_i$ are orthogonal, the above matrices
are pairs of orthogonal matrices and so are orthogonal.
Note that each of the $2^{2^n-n-1}$ $A_i's$ is paired with all the other $2^{2^n-n-1}$
$A_j's$ and each is paired with $2^n$ shifts of the rows.  So the total number of matrices above is:
\[ 2^{2^n-n-1}\cdot 2^{2^n-n-1}\cdot 2^n= 
= 2^{2^{n+1}-(n+1)-1}.
\]
I.e.  We have used up all the distinct Hadamard vectors
in a space of dimension $2^{n+1}$.

\end{document}